# Empirical and Gaussian processes on Besov classes

**Richard Nickl**

*Department of Statistics, University of Vienna*

**Abstract:** We give several conditions for pregaussianity of norm balls of Besov spaces defined over $\mathbb{R}^d$ by exploiting results in Haroske and Triebel (2005). Furthermore, complementing sufficient conditions in Nickl and Pötscher (2005), we give necessary conditions on the parameters of the Besov space to obtain the Donsker property of such balls. For certain parameter combinations Besov balls are shown to be pregaussian but not Donsker.

## 1. Introduction

Bounds for the size (measured, e.g., by metric entropy) of a subset $\mathcal{F}$ of the space $\mathcal{L}^2(\mathbb{P})$ of functions square-integrable w.r.t. some probability measure $\mathbb{P}$ allow one to derive limit theorems for the empirical process (indexed by $\mathcal{F}$) as well as continuity properties of the (limiting) Gaussian process (indexed by $\mathcal{F}$). These bounds are often derived from smoothness conditions on the functions contained in $\mathcal{F}$. Function classes that satisfy differentiability or Hölder conditions were among the first examples for pregaussian and Donsker classes, cf. Strassen and Dudley [14], Giné [7], Stute [15], Marcus [11], Giné and Zinn [8], Arcones [1] and van der Vaart [18]. In recent years, interest in spaces of functions with 'generalized smoothness', e.g., spaces of Besov- and Triebel- type, has grown. These spaces contain the spaces defined by more classical smoothness conditions (such as Hölder(-Zygmund), Lipschitz and Sobolev spaces) as special cases and serve as a unified theoretical framework. Besov and Triebel spaces play an increasing role in nonparametric statistics, information theory and data compression, see, e.g., Donoho and Johnstone [3], Donoho, Vetterli, DeVore and Daubechies [4] and Birgé and Massart [2]. Relatively little was known until recently about empirical and Gaussian processes on such function classes, in particular with focus on spaces defined over the whole Euclidean space $\mathbb{R}^d$. Building on Haroske and Triebel [10], sufficient conditions for the parameters of the Besov space were given in [12] implying that the corresponding norm balls are Donsker classes. In the present paper, we extend and complement these results. We give necessary and sufficient conditions for the pregaussian/Donsker property of balls in Besov spaces. In certain 'critical' cases, Besov balls are shown to be pregaussian but *not* Donsker.

## 2. Besov spaces

For $h$ a real-valued Borel-measurable function defined on $\mathbb{R}^d$ ($d \in \mathbb{N}$) and $\mu$ a (nonnegative) Borel measure on $\mathbb{R}^d$, we set $\mu f := \int_{\mathbb{R}^d} f d\mu$ as well as $\|h\|_{r,\mu} := (\int_{\mathbb{R}^d} |h|^r d\mu)^{1/r}$ for $1 \leq r \leq \infty$ (where $\|h\|_{\infty,\mu}$ denotes the $\mu$-essential supremum of $|h|$). As usual, we denote by $\mathcal{L}^r(\mathbb{R}^d, \mu)$ the vector space of all Borel-measurable functions $h : \mathbb{R}^d \to \mathbb{R}$ that satisfy $\|h\|_{r,\mu} < \infty$. In accordance, $L^r(\mathbb{R}^d, \mu)$ denotes the corresponding Banach spaces of equivalence classes $[h]_\mu$, $h \in \mathcal{L}^r(\mathbb{R}^d, \mu)$, modulo equality $\mu$-a.e. The symbol $\lambda$ will be used to denote Lebesgue-measure on $\mathbb{R}^d$.

We follow Edmunds and Triebel ([6], 2.2.1) in defining Besov spaces: Let $\varphi_0$ be a complex-valued $C^\infty$-function on $\mathbb{R}^d$ with $\varphi_0(x) = 1$ if $\|x\| \leq 1$ and $\varphi_0(x) = 0$ if $\|x\| \geq 3/2$. Define $\varphi_1(x) = \varphi_0(x/2) - \varphi_0(x)$ and $\varphi_k(x) = \varphi_1(2^{-k+1}x)$ for $k \in \mathbb{N}$. Then the functions $\varphi_k$ form a dyadic resolution of unity. Let $\mathcal{S}(\mathbb{R}^d)$ denote the Schwartz space of rapidly decreasing infinitely differentiable complex-valued functions and let $\mathcal{S}'(\mathbb{R}^d)$ denote the (dual) space of complex tempered distributions on $\mathbb{R}^d$.. In this paper we shall restrict attention to real-valued tempered distributions $T$ (i.e., $T = \bar{T}$, where $\bar{T}$ is defined via $\bar{T}(\phi) = \overline{T(\bar{\phi})}$ for $\phi \in \mathcal{S}(\mathbb{R}^d)$). Let $F$ denote the Fourier transform acting on $\mathcal{S}'(\mathbb{R}^d)$ (see, e.g., Chapter 7.6 in [13]). Then $F^{-1}(\varphi_k FT)$ is an entire analytic function on $\mathbb{R}^d$ for any $T \in \mathcal{S}'(\mathbb{R}^d)$ and any $k$ by the Paley-Wiener-Schwartz theorem (see, e.g., p. 272 in [13]).

**Definition 1** (Besov spaces). Let $-\infty < s < \infty$, $1 \leq p \leq \infty$, and $1 \leq q \leq \infty$. For $T \in \mathcal{S}'(\mathbb{R}^d)$ define

$$\|T\|_{s,p,q,\lambda} := \left( \sum_{k=0}^\infty 2^{ksq} \left\| F^{-1}(\varphi_k FT) \right\|_{p,\lambda}^q \right)^{1/q}$$

with the modification in case $q = \infty$

$$\|T\|_{s,p,\infty,\lambda} := \sup_{0 \leq k < \infty} 2^{ks} \left\| F^{-1}(\varphi_k FT) \right\|_{p,\lambda}.$$

Define further (real) Besov spaces as

$$B^s_{pq}(\mathbb{R}^d) := \{T \in \mathcal{S}'(\mathbb{R}^d) : T = \bar{T}, \ \|T\|_{s,p,q,\lambda} < \infty\}.$$

$B^s_{pq}(\mathbb{R}^d)$ is a Banach space and the norm is independent of the choice of $\varphi_0$, and, in particular, different $\varphi_0$ result in equivalent norms, cf. Edmunds and Triebel [6], 2.2.1.

**Remark 2.** (i) The focus in the present paper will be on $s > 0$, in which case it follows (e.g., from 2.3.2 in [17]) that $B^s_{pq}(\mathbb{R}^d)$ consists of (equivalence classes of) $p$-fold integrable functions. In fact, for these parameters, we could alternatively have defined the spaces $B^s_{pq}(\mathbb{R}^d)$ as $\{[f]_\lambda \in L^p(\mathbb{R}^d, \lambda), \|f\|_{s,p,q,\lambda} < \infty\}$.

(ii) We note that $\|T\|_{s,p,q,\lambda} < \infty$ if and only if $\|\bar{T}\|_{s,p,q,\lambda} < \infty$ for any $T \in \mathcal{S}'(\mathbb{R}^d)$. In fact,

$$\|T\|_{s,p,q,\lambda} \leq \|\mathrm{Re}T\|_{s,p,q,\lambda} + \|\mathrm{Im}T\|_{s,p,q,\lambda} \leq c \|T\|_{s,p,q,\lambda}$$

holds for some $1 \leq c < \infty$ and for every $T \in \mathcal{S}'(\mathbb{R}^d)$. As a consequence, one can easily carry over results for complex Besov spaces to real ones and vice versa.

(iii) At least for positive $s$, there are many equivalent norms on $B^s_{pq}(\mathbb{R}^d)$, some of them possibly more common than the one used in Definition 1; see, e.g., Remark 2 in [12]. In particular, the Hölder-Zygmund Spaces are identical (up to an equivalent norm) to the spaces $B^s_{\infty\infty}(\mathbb{R}^d)$ if $s > 0$.

(iv) Triebel spaces $F^s_{pq}(\mathbb{R}^d)$ are defined in 2.2.1/7 in [6]. We have the chain of continuous imbeddings $B^s_{pu}(\mathbb{R}^d) \hookrightarrow F^s_{pq}(\mathbb{R}^d) \hookrightarrow B^s_{pv}(\mathbb{R}^d)$ for $0 < u \leq \min(p,q)$ and $\max(p,q) \leq v \leq \infty$. By using these imbeddings, the results of the present paper can also be applied to Triebel spaces. Note that $F_{p2}^0(\mathbb{R}^d) = L^p(\mathbb{R}^d, \lambda)$ holds, and for positive $s$, we have that $F_{p2}^s(\mathbb{R}^d)$ is equal to the classical Sobolev spaces. See 2.2.2 in [6] for further details.

Let $\mathsf{C}(\mathbb{R}^d)$ be the vector space of bounded continuous real-valued functions on $\mathbb{R}^d$ normed by the sup-norm $\|\cdot\|_\infty$. If either $s > d/p$ or $s = d/p$ and $q = 1$, it is well-known (see, e.g., Proposition 3 in [12]) that each equivalence class $[f]_\lambda \in B_{pq}^s(\mathbb{R}^d)$, contains a (unique) continuous representative. [In fact, the Banach space $B_{pq}^s(\mathbb{R}^d)$ is imbedded (up to a section map) into the space $\mathsf{C}(\mathbb{R}^d)$.] Hence, if either $s > d/p$ or $s = d/p$ and $q = 1$, we can define the (closely related) Banach space

$$\mathsf{B}_{pq}^s(\mathbb{R}^d) = \{f \in \mathsf{C}(\mathbb{R}^d) : [f]_\lambda \in L^p(\mathbb{R}^d, \lambda), \|f\|_{s,p,q,\lambda} < \infty\}$$

(again normed by $\|\cdot\|_{s,p,q,\lambda}$) by collecting the continuous representatives.

Throughout the paper we shall use the following notational agreements: We define the function $\langle x \rangle^\gamma = (1 + \|x\|^2)^{\gamma/2}$ parameterized by $\gamma \in \mathbb{R}$, where $x$ is an element of $\mathbb{R}^d$ and where $\|\cdot\|$ denotes the Euclidean norm. Also, for two real-valued functions $a(\cdot)$ and $b(\cdot)$, we write $a(\varepsilon) \lesssim b(\varepsilon)$ if there exists a positive (finite) constant $c$ not depending on $\varepsilon$ such that $a(\varepsilon) \leq cb(\varepsilon)$ holds for all $\varepsilon > 0$. If $a(\varepsilon) \lesssim b(\varepsilon)$ and $b(\varepsilon) \lesssim a(\varepsilon)$ both hold we write $a(\varepsilon) \sim b(\varepsilon)$. [In abuse of notation, we shall also use this notation for sequences $a_k$ and $b_k$, $k \in \mathbb{N}$ as well as for two (semi)norms $\|\cdot\|_{X,1}$ and $\|\cdot\|_{X,2}$ on a vector space $X$.]

## 3. Main results

Let $(S, \mathcal{A}, \mu)$ be some probability space and let $\mathbb{P}$ be a (Borel) probability measure on $\mathbb{R}^d$. Let $\emptyset \neq \mathcal{F} \subseteq \mathcal{L}^2(\mathbb{R}^d, \mathbb{P})$. A Gaussian process $\mathbb{G} : (S, \mathcal{A}, \mu) \times \mathcal{F} \to \mathbb{R}$ with mean zero and covariance $E\mathbb{G}(f)\mathbb{G}(g) = \mathbb{P}[(f - \mathbb{P}f)(g - \mathbb{P}g)]$ for $f, g \in \mathcal{F}$ is called a (generalized) *Brownian bridge* process on $\mathcal{F}$. The covariance induces a semimetric $\rho^2(f, g) = E[\mathbb{G}(f) - \mathbb{G}(g)]^2$ for $f, g \in \mathcal{F}$. A function class $\mathcal{F} \subseteq \mathcal{L}^2(\mathbb{R}^d, \mathbb{P})$ will be called $\mathbb{P}$-*pregaussian* if such a Gaussian process $\mathbb{G}$ can be defined such that for every $s \in S$, the map $f \longmapsto \mathbb{G}(f, s)$ is bounded and uniformly continuous w.r.t. the semimetric $\rho$ from $\mathcal{F}$ into $\mathbb{R}$. For further details see p.92-93 in [5].

Let $\mathbb{P}_n = 1/n \sum_{i=1}^n \delta_{X_i}$ denote the empirical measure of $n$ independent $\mathbb{R}^d$-valued random variables $X_1, \ldots, X_n$ identically distributed according to some law $\mathbb{P}$. [We assume here the standard (canonical) model as on p.91 in [5].] For $\mathcal{F} \subseteq \mathcal{L}^2(\mathbb{R}^d, \mathbb{P})$, the $\mathcal{F}$-indexed *empirical process* $\nu_n$ is given by $f \longmapsto \nu_n(f) = \sqrt{n}\,(\mathbb{P}_n - \mathbb{P})\,f$. The class $\mathcal{F}$ is said to be $\mathbb{P}$-*Donsker* if it is $\mathbb{P}$-pregaussian and if $\nu_n$ converges in law in the space $\ell^\infty(\mathcal{F})$ to a (generalized) Brownian bridge process over $\mathcal{F}$, cf. p.94 in [5]. Here $\ell^\infty(\mathcal{F})$ denotes the Banach space of all bounded real-valued functions on $\mathcal{F}$. If $\mathcal{F}$ is $\mathbb{P}$-Donsker for all probability measures $\mathbb{P}$ on $\mathbb{R}^d$, it is called *universally Donsker*.

In [12], Corollary 5, Proposition 1 and Theorem 2, the following results were proved. [Clearly, one may replace $\mathcal{U}$ by and bounded subset of $\mathsf{B}_{pq}^s(\mathbb{R}^d)$ in the proposition.]

**Proposition 3.** *Let $\mathcal{U}$ be the closed unit ball of $\mathsf{B}_{pq}^s(\mathbb{R}^d)$ where $1 \leq p \leq \infty$, $1 \leq q \leq \infty$. Let $\mathbb{P}$ be a probability measure on $\mathbb{R}^d$.*

1. *Let $1 \leq p \leq 2$ and $s > d/p$. Then $\mathcal{U}$ is $\mathbb{P}$-Donsker, and hence also $\mathbb{P}$-pregaussian.*

2. Let $2 < p \leq \infty$ and $s > d/2$. Assume that $\int_{\mathbb{R}^d} \|x\|^{2\gamma} d\mathbb{P} < \infty$ holds for some $\gamma > d/2 - d/p$. Then $\mathcal{U}$ is $\mathbb{P}$-Donsker, and hence also $\mathbb{P}$-pregaussian.
3. Let $d = q = 1$, $1 \leq p < 2$ and $s = 1/p$. Then $\mathcal{U}$ is $\mathbb{P}$-Donsker, and hence also $\mathbb{P}$-pregaussian.

In the present paper we show on the one hand that, if one is interested in the pregaussian property only, the conditions of Proposition 3 can be substantially weakened. On the other hand, we show that Proposition 3 is (essentially) best possible w.r.t. the Donsker property: It turns out that $s \geq \max(d/p, d/2)$ always has to be satisfied for $\mathcal{U}$ to be $\mathbb{P}$-Donsker and that the moment condition in Part 2 of Proposition 3 cannot be improved upon. We also give a rather definite picture of the limiting case $s = d/p$ (where only the cases $q = 1$ and $d > 1$, as well as $p = 2$ and $q = 1$, will remain undecided).

### 3.1. The pregaussian property

We first discuss the pregaussian property in the 'nice' case $s > \max(d/p, d/2)$: If $s > d/p$ and $p \leq 2$, Proposition 3 implies that the unit ball of the Besov space is pregaussian for every probability measure. On the other hand, maybe not surprisingly, if the integrability parameter $p$ of the Besov space is larger than 2, Proposition 1 requires an additional moment condition on the probability measure to obtain the pregaussian property. The following theorem shows that this additional moment condition is also necessary (for most probability measures possessing Lebesgue-densities). [Note that $s > d/2$ ensures also that $s > d/p$ holds, so the condition $s > \max(d/p, d/2)$ is always satisfied.]

**Theorem 4.** *Let $\mathcal{U}$ be the closed unit ball of $\mathsf{B}_{pq}^s(\mathbb{R}^d)$ with $2 < p \leq \infty$, $1 \leq q \leq \infty$ and $s > d/2$. Let $\delta$ be arbitrary subject to $0 < \delta \leq d/2 - d/p$. Define the probability measure $\mathbb{P}$ by $d\mathbb{P}(x) = \varphi(x) \langle x \rangle^{-d-2\delta} d\lambda(x)$ where $0 < c \leq \varphi(x)$ holds for some constant $c$ and all $x \in \mathbb{R}^d$ (and where $\|\varphi \langle x \rangle^{-d-\delta}\|_{1,\lambda} = 1$). Then the set $\mathcal{U}$ is not $\mathbb{P}$-pregaussian.*

*Proof.* Note first that $\mathcal{U}$ is a bounded subset of $\mathsf{C}(\mathbb{R}^d)$ (see, e.g., Proposition 3 in [12]) and hence also of $\mathcal{L}^r(\mathbb{R}^d, \mathbb{P})$ for every $1 \leq r \leq \infty$. Observe that

$$\|f - g\|_{2,\mathbb{P}}^2 = \int_{\mathbb{R}^d} [f - g]^2 d\mathbb{P} = \int_{\mathbb{R}^d} [(f - g) \langle x \rangle^{-(d-2\delta)/2}]^2 \varphi d\lambda$$
$$\geq c \left\| (f - g) \langle x \rangle^{-(d-2\delta)/2} \right\|_{2,\lambda}^2$$

holds for $f, g \in \mathcal{L}^2(\mathbb{R}^d, \mathbb{P})$. Hence we have for the metric entropy (see Definition 9 in the Appendix) that

$$H(\varepsilon, U, \|\cdot\|_{2,\mathbb{P}}) \geq H(\varepsilon/c, U, \left\| (\cdot) \langle x \rangle^{-(d-2\delta)/2} \right\|_{2,\lambda})$$

holds. We obtain a lower bound of order $\varepsilon^{-\alpha}$ for the r.h.s. of the above display from Corollary 12 in the Appendix upon setting $\gamma = (d - 2\delta)/2$ in that corollary. Since $s - d/p > d/2 - d/p > 0$ and $\delta \leq d/2 - d/p$, it follows that $\gamma < s - d/p + d/2$ and we obtain $\alpha = (\delta/d + 1/p)^{-1}$. Clearly $\alpha > 2$ holds since $\delta \leq d/2 - d/p$. Define the Gaussian process $\mathbb{L}(f) = \mathbb{G}(f) + Z \cdot \mathbb{P}f$ for $f \in \mathcal{U}$ where $Z$ is a standard normal variable independent of $\mathbb{G}$. It is easily seen that this process has covariance $E\mathbb{L}(f)\mathbb{L}(g) = \mathbb{P}fg$. Since $a > 2$ and since $\mathbb{P}$ possesses a Lebesgue-density, we can apply the Sudakov-Chevet minoration (Theorem 2.3.5 in [5]) which implies that the process $\mathbb{L}$ is $\mu$-a.s. unbounded on $\mathcal{U}$. Since $\sup_{f \in \mathcal{U}} |\mathbb{P}f| < \infty$ holds, we have that $\sup_{f \in \mathcal{U}} |\mathbb{L}(f)| = \infty$ $\mu$-a.s. implies $\sup_{f \in \mathcal{U}} |\mathbb{G}(f)| = \infty$ $\mu$-a.s. This proves that $\mathcal{U}$ is not $\mathbb{P}$-pregaussian. □

The set $\mathcal{U}$ is uniformly bounded (in fact, for $p < \infty$, any $f \in \mathcal{U}$ satisfies $\lim_{\|x\| \to \infty} f(x) = 0$, see Proposition 3 in [12]), but nevertheless one needs a moment condition on the probability measure to obtain the pregaussian property. [The reason is, not surprisingly, that the degree of compactness in $L^2(\mathbb{R}^d, \mathbb{P})$ measured in terms of metric entropy is driven *both* by smoothness of the function class *and* by its rate of decay at infinity.]

In the remainder of this section we shed light on the critical cases $s \leq d/p$ and/or $s \leq d/2$. The following proposition shows that in case $s \leq d/p$ but $s > d/2$ (and hence $1 \leq p < 2$), Besov balls are again pregaussian for a large class of probability measures:

**Theorem 5.** *Let $U$ be the closed unit ball of $B_{pq}^s(\mathbb{R}^d)$ with $1 \leq p < 2$, $1 \leq q \leq \infty$ and $s > d/2$, and let $\mathcal{U}$ be any set constructed by selection of one arbitrary representative out of every $[f]_\lambda \in U$. Let $\mathbb{P}$ be a probability measure on $\mathbb{R}^d$ that possesses a density $\varphi$ w.r.t. Lebesgue measure on $\mathbb{R}^d$ such that $\|\varphi \langle x \rangle^d\|_{\infty, \lambda} < \infty$. Then $\mathcal{U}$ is $\mathbb{P}$-pregaussian.*

*Proof.* Note first that $U$ is a bounded subset of $L^2(\mathbb{R}^d, \lambda)$ (by Proposition 11 and (4) in the Appendix) and hence also of $L^2(\mathbb{R}^d, \mathbb{P})$ since $[\varphi]_\lambda \in L^\infty(\mathbb{R}^d, \lambda)$. Observe next that

$$\|f - g\|_{2,\mathbb{P}}^2 = \int_{\mathbb{R}^d} [f - g]^2 \varphi d\lambda = \int_{\mathbb{R}^d} [(f - g) \langle x \rangle^{-d/2}]^2 \varphi \langle x \rangle^d \, d\lambda$$
$$\leq \left\|(f - g) \langle x \rangle^{-d/2}\right\|_{2,\lambda}^2 \left\|\varphi \langle x \rangle^d\right\|_{\infty, \lambda}$$

holds for $f, g \in \mathcal{L}^2(\mathbb{R}^d, \mathbb{P})$ by Hölder's inequality. Hence we apply Corollary 12 in the Appendix with $\gamma = d/2$ to obtain

$$H(\varepsilon, \mathcal{U} \, \|\cdot\|_{2,\mathbb{P}}) \leq H(\varepsilon \left\|\varphi \langle x \rangle^d\right\|_{\infty, \lambda}, U, \left\|(\cdot) \langle x \rangle^{-d/2}\right\|_{2,\lambda}) \lesssim \varepsilon^{-\alpha}$$

where $\alpha = d/s$ if $s - d/p < 0$ and $\alpha = p$ if $s - d/p > 0$ and where we have used that $\mathbb{P}$ is absolutely continuous w.r.t. Lebesgue measure $\lambda$. In both cases we have $\alpha < 2$. Hence we can apply Theorem 2.6.1 in [5] to obtain (a.s.) sample-boundedness and -continuity of the process $\mathbb{L}$ (defined in the proof of Theorem 4 above) on $\mathcal{U}$ w.r.t. the $\mathcal{L}^2(\mathbb{R}^d, \mathbb{P})$-seminorm. If $\pi(f) = f - \mathbb{P}f$, then $\mathbb{L}(\pi(f)) = \mathbb{G}(f)$ is also (a.s.) sample-bounded and -continuous on $\mathcal{U}$ and hence we obtain the $\mathbb{P}$-pregaussian property for $\mathcal{U}$ by the same reasoning as on p.93 in [5]. If $s = d/p$, view $U$ as a bounded subset of $B_{pq}^{d/p - \varepsilon}(\mathbb{R}^d)$ where $\varepsilon$ can be chosen small enough such that $d/p - \varepsilon > d/2$ holds (note that $d/p > d/2$ since $p < 2$) and hence the pregaussian property follows from the case $s - d/p < 0$ just established. This finishes the proof. □

Note that any probability measure $\mathbb{P}$ that possesses a bounded density which is eventually monotone, or a bounded density with polynomial or exponential tails, satisfies the condition of the theorem. [We note that at least for the special case $d = q = 1$, $s = 1/p$, $p < 2$, the condition on $\mathbb{P}$ can be removed by Proposition 3.]

The following theorem deals with the remaining cases and shows that $s \geq d/2$ always has to be satisfied (irrespective of $p$) to obtain the pregaussian property.

**Theorem 6.** *Let $U$ be the closed unit ball of $B^s_{pq}(\mathbb{R}^d)$ with $1 \le p \le \infty$, $1 \le q \le \infty$, $0 < s < d/2$, and let $\mathcal{U}$ be any set constructed by selection of one arbitrary representative out of every $[f]_\lambda \in U$. Let $\mathbb{P}$ be a probability measure that possesses a bounded density $\varphi$ w.r.t. Lebesgue measure on $\mathbb{R}^d$ which satisfies $0 < c \le \varphi(x)$ for some constant $c$ and all $x$ in some open subset $V$ of $\mathbb{R}^d$. Then the set $\mathcal{U}$ is not $\mathbb{P}$-pregaussian.*

*Proof.* Since $V$ is open, it contains an open Euclidean ball $\Omega$, which is a bounded $C^\infty$-domain in the sense of Triebel [17], 3.2.1. Denote by $\lambda\,|\Omega$ Lebesgue measure on $\Omega$ and by $L^2(\Omega, \lambda)$ the usual Banach space normed by the usual $L^2$-norm $\|\cdot\|_{2,\lambda|\Omega}$ on $\Omega$. Let $U\,|\Omega$ be the set of restrictions $[f\,|\Omega]_{\lambda|\Omega}$ of elements $[f]_\lambda \in U$ to the set $\Omega$. Note that $U\,|\Omega$ is the unit ball of the factor Besov space $B^s_{pq}(\mathbb{R}^d)\,|\Omega$ over $\Omega$ obtained by restricting the elements of $B^s_{pq}(\mathbb{R}^d)$ to $\Omega$ with the restricted Besov norm

$$\|f\|_{s,p,q,|\Omega} := \inf\left\{\|g\|_{s,p,q,\lambda}:\; [g]_\lambda \in B^s_{pq}(\mathbb{R}^d),\quad [g\,|\Omega]_{\lambda|\Omega} = f\right\}.$$

We first handle the case $p = 1$. In view of 2.5.1/7 and 2.2.2/1 in [6], we have that $B^{d/2}_{1\infty}(\mathbb{R}^d)\,|\Omega \not\subseteq L^2(\Omega, \lambda)$. But by $s < d/2$ and 3.3.1/7 of Triebel (1983) we also have $B^{d/2}_{1\infty}(\mathbb{R}^d)\,|\Omega \subseteq B^s_{1q}(\mathbb{R}^d)\,|\Omega$ hence we conclude that $B^s_{1q}(\mathbb{R}^d)\,|\Omega \not\subseteq L^2(\Omega, \lambda)$. Since $\varphi \ge c$ on $\Omega$, this implies $U \not\subseteq L^2(\mathbb{R}^d, \mathbb{P})$, so $\mathcal{U}$ cannot be $\mathbb{P}$-pregaussian.

We now turn to $p > 1$. We first treat the case $s = d/2 - \varepsilon$ where $\varepsilon > 0$ is arbitrary subject to $\varepsilon < d - d/p$. Then $U$ is a bounded subset of $L^2(\mathbb{R}^d, \lambda)$ by Proposition 11 and (4) in the Appendix and hence also of $L^2(\mathbb{R}^d, \mathbb{P})$ since $\varphi$ is bounded. We now obtain a metric entropy lower bound for $U$ in $L^2(\mathbb{R}^d, \mathbb{P})$. Observe that

$$\|f - g\|^2_{2,\mathbb{P}} = \int_{\mathbb{R}^d} [f - g]^2 d\mathbb{P} \ge c \int_\Omega [f - g]^2 d\lambda\,|\Omega$$

holds for $f, g \in \mathcal{L}^2(\mathbb{R}^d, \mathbb{P})$ and hence

(1) $$H(\varepsilon, U, \|\cdot\|_{2,\mathbb{P}}) \ge H(\varepsilon/c, U\,|\Omega, \|\cdot\|_{2,\lambda|\Omega})$$

holds. By 3.3.3/1 in [6], we obtain the entropy number (see Definition 8 in the Appendix)

$$e\left(k, id(U\,|\Omega), \|\cdot\|_{0,2,\infty,|\Omega}\right) \sim k^{-s/d}.$$

Now by Lemma 1 as well as expression (4) in the Appendix we obtain

$$H(\varepsilon, U\,|\Omega, \|\cdot\|_{2,\lambda|\Omega}) \gtrsim \varepsilon^{-d/s}.$$

But since $s < d/2$ holds by assumption, this (together with (1)) implies that $\sup_{f \in \mathcal{U}} |\mathbb{G}(f)| = \infty$ $\mu$-a.s. by the same application of the Sudakov-Chevet minoration as in the proof of Theorem 4 above, noting that $\sup_{f \in \mathcal{U}} |\mathbb{P}f| < \infty$ holds since $U$ is bounded in $L^2(\mathbb{R}^d, \mathbb{P})$. Hence $\mathcal{U}$ is not $\mathbb{P}$-pregaussian in this case. The remaining cases $s - \varepsilon$ with $\varepsilon \ge d - d/p$ now follow from the continuous imbedding $B^s_{pq}(\mathbb{R}^d) \hookrightarrow B^t_{pq}(\mathbb{R}^d)$ for $s > t$, cf. 2.3.2/7 in [17]. □

Observe that $\mathbb{P}$ in the above theorem could be compactly supported, so the pregaussian property cannot be restored by a moment condition. Inspection of the proof shows that a similar negative result can be proved for the unit ball of a Besov space over any subdomain of $\mathbb{R}^d$ (that possesses a suitably regular boundary).

The limiting case $p = 2$, $1 \leq q \leq \infty$, $s = d/2$ remains open: Here, one would have to go to the logarithmic scale of metric entropy rates, in which case it is known that metric entropy conditions are not sharp in terms of proving the pregaussian property, see p.54 in [5]. At least for $q \geq 2$ we conjecture that the unit ball of $B_{2q}^{d/2}(\mathbb{R}^d)$ is not $\mathbb{P}$-pregaussian for absolutely continuous probability measures possessing a bounded density.

## 3.2. The Donsker property

In this section we show that Proposition 3 is (essentially) best possible in terms of the Donsker property for norm balls in Besov spaces. We first discuss the 'nice' case $s > \max(d/p, d/2)$. If $p \leq 2$, Part of Proposition 3 is certainly best possible (since then Besov balls are universally Donsker). Since $\mathbb{P}$-Donsker classes must be $\mathbb{P}$-pregaussian, the moment condition in Part 2 ($p > 2$) of Proposition 3 is (essentially) necessary in view of Theorem 4 above. For the case $p = q = \infty$, these findings imply known results for Hölder and Lipschitz classes due to Giné and Zinn [8], Arcones [1] and van der Vaart [18]; cf. also the discussion in Remark 5 in [12].

We now turn to the critical cases $s \leq d/p$ and/or $s \leq d/2$. Since Donsker classes need to be pregaussian, Theorem 6 implies that $s \geq d/2$ *always has to be satisfied* (at least for the class of probability measures defined in that theorem).

On the other hand, for $1 \leq p < 2$ we showed in Theorem 5 that norm balls of $B_{pq}^s(\mathbb{R}^d)$ with $d/2 < s \leq d/p$ (and hence $1 \leq p < 2$) are pregaussian for a large class of probability measures. So the question arises whether these classes are also Donsker classes for these probability measures. In the special case $q = d = 1$ and $s = 1/p$, these classes are in fact universally Donsker in view of Part 3 of Proposition 3. We do not know whether this can be generalized to the case $d > 1$, that is, whether the unit ball of $\mathsf{B}_{p1}^{d/p}(\mathbb{R}^d)$ with $1 \leq p < 2$ is a (universal) Donsker class. [The proof in case $d = 1$ uses spaces of functions of bounded $p$-variation, a concept which is not straightforwardly available for $d > 1$.]

On the other hand, the following theorem shows that the function classes that were shown to be pregaussian in Theorem 5 are in fact *not* $\mathbb{P}$-Donsker for probability measures $\mathbb{P}$ possessing a bounded density if $s < d/p$, or if $s = d/p$ but $q > 1$ hold. The proof strategy partially follows the proof of Theorem 2.3 in [11].

**Theorem 7.** *Let $U$ be the closed unit ball of $B_{pq}^s(\mathbb{R}^d)$ with $1 \leq p \leq \infty$, $1 \leq q \leq \infty$, $s > 0$ and let $\mathcal{U}$ be any set constructed by selection of one arbitrary representative out of every $[f]_\lambda \in U$. Assume that $\mathbb{P}$ possesses a bounded density w.r.t. Lebesgue measure. Suppose that either $s < d/p$ or that $s = d/p$ but $q > 1$ holds. Then $\mathcal{U}$ is not a $\mathbb{P}$-Donsker class.*

*Proof.* We first consider the case $s = d/p$ but $q > 1$. By Theorem 2.6.2/1 in [16], p. 135, $B_{pq}^{d/p}(\mathbb{R}^d)$ contains a function $\psi \in \mathcal{L}^1(\mathbb{R}^d, \lambda)$ that satisfies

$$|\psi(x)| \geq C \log |\log |x||$$

for $|x| \in (0, \varepsilon]$ and some $0 < \varepsilon < 1$. We may assume w.l.o.g. $\|\psi\|_{s,p,q,\lambda} \leq 1$. Since $(F\psi(\cdot - y))(u) = e^{-iyu} F\psi(u)$ holds, inspection of Definition 1 shows that $\|\psi(\cdot - y)\|_{s,p,q,\lambda} = \|\psi\|_{s,p,q,\lambda} \leq 1$ for every $y \in \mathbb{R}^d$. Let $(z_i)_{i=1}^\infty$ denote all points in $\mathbb{R}^d$ with rational coordinates and define $\psi_i = \psi(\cdot - z_i)$ which satisfies $\|\psi_i\|_{s,p,q,\lambda} \leq 1$ for every $i$. Consequently, we have $\{\tilde{\psi}_i\}_{i=1}^\infty \subseteq \mathcal{U}$ where $\tilde{\psi}_i$ is obtained by modifying each $\psi_i$ on a set $N_i$ of Lebesgue-measure zero if necessary. Clearly $\cup_{i=1}^\infty N_i$ is again a set of Lebesgue measure zero. Let now $x \in \mathbb{R}^d \setminus \cup_{i=1}^\infty N_i$ be arbitrary and let the index set $I_x$ consist of all $i \in \mathbb{N}$ s.t. $|x - z_i| < \varepsilon$ holds. Clearly $(z_i)_{i \in I_x}$ is dense in a neighborhood of $x$. Consequently

$$\sup_{f \in \mathcal{U}} |f(x)| \geq \sup_{i \in \mathbb{N}} \left|\tilde{\psi}_i(x)\right| = \sup_{i \in \mathbb{N}} |\psi_i(x)| = \sup_{i \in \mathbb{N}} |\psi(x - z_i)|$$

$$\geq \sup_{i \in I_x} C \log |\log |x - z_i|| = \infty$$

holds for every $x \in \mathbb{R}^d \setminus \cup_{i=1}^\infty N_i$ and hence Lebesgue almost everywhere. Note furthermore that $U$ is bounded in $L^2(\mathbb{R}^d, \lambda)$ (by Proposition 11 and (4) in the Appendix). Furthermore, $\mathbb{P}$ possesses a density $[\phi]_\lambda \in L^\infty(\mathbb{R}^d, \lambda) \cap L^1(\mathbb{R}^d, \lambda) \subseteq L^2(\mathbb{R}^d, \lambda)$, so we have $\sup_{f \in \mathcal{U}} |\mathbb{P}f| < \infty$ by using the Cauchy-Schwarz inequality. Conclude that

$$M_\mathbb{P}(x) = \sup_{f \in \mathcal{U}} |f(x) - \mathbb{P}f| \geq \sup_{f \in \mathcal{U}} |f(x)| - \sup_{f \in \mathcal{U}} |\mathbb{P}f| = \infty$$

holds $\lambda$-a.e. Since $\mathbb{P}$ is absolutely continuous, we have that $\mathcal{U}$ is not a $\mathbb{P}$-Donsker class since $t^2 \mathbb{P}(M_\mathbb{P} > t) \to 0$ is necessary for the $\mathbb{P}$-Donsker property to hold for $\mathcal{U}$ (see, e.g., Proposition 2.7 in [9]). The remaining cases follow from the continuous imbedding $B_{pu}^{d/p}(\mathbb{R}^d) \hookrightarrow B_{pv}^s(\mathbb{R}^d)$ for $s < d/p$ and $u, v \in [1, \infty]$ (cf. 2.3.2/7 in [17]). □

At least on the sample space $\mathbb{R}^d$ we are not aware of any other ('constructive') examples for pregaussian classes that are not Donsker: The above theorem shows that the empirical process does *not* converge in law in $\ell^\infty(\mathcal{U})$ if $\mathcal{U}$ is the unit ball of $B_{pq}^s(\mathbb{R}^d)$ (with $s < d/p$ or $s = d/p$ but $q > 1$). However, if $p < 2$ and $s > d/2$ a sample-bounded and -continuous Brownian bridge process *can* be defined on $\mathcal{U}$ by Theorem 5 above.

Inspection of the proof shows that a similar negative result can be proved for the unit ball of a Besov space defined over any (non-empty) subset $\Omega$ of $\mathbb{R}^d$ (at least if $\Omega$ has regular boundary). Note that the above theorem also implies for the case $p = 2$, $s = d/2$ (not covered in Section 3.1) that the unit ball of $B_{2q}^{d/2}(\mathbb{R}^d)$ is not Donsker if $q > 1$. [The special case $q = 1$ remains open.]

## Appendix A: Technical results

**Definition 8.** Let $\mathcal{J}$ be a subset of the normed space $(Y, \|\cdot\|_Y)$, and let $U_Y = \{y \in Y : \|y\|_Y \leq 1\}$ be the closed unit ball in $Y$. Then, for all natural numbers $k$, the k-th entropy number of $\mathcal{J}$ is defined as

$$e(k, \mathcal{J}, \|\cdot\|_Y) = \inf \left\{ \varepsilon > 0 : \mathcal{J} \subseteq \bigcup_{j=1}^{2^{k-1}} (y_j + \varepsilon U_Y) \text{ for some } y_1, \ldots, y_{2^{k-1}} \in Y \right\},$$

with the convention that the infimum equals $+\infty$ if the set over which it is taken is empty.

Suppose $(X, \|\cdot\|_X)$ and $(Y, \|\cdot\|_Y)$ are normed spaces such that $X$ is a linear subspace of $Y$. Let $U_X$ the closed unit ball in $X$. Then, $e(k, id(U_X), \|\cdot\|_Y)$ is called the k-th entropy number of the operator $id : X \to Y$. Clearly, $e(k, id(U_X), \|\cdot\|_Y)$ is finite for all $k$ if and only if $id$ is continuous from $X$ to $Y$ (in which case we shall write $(X, \|\cdot\|_X) \hookrightarrow (Y, \|\cdot\|_Y)$) and the entropy numbers converge to zero as $k \to \infty$ if and only if the operator $id$ is compact (has totally bounded image in $Y$.)

**Definition 9.** For a (non-empty) subset $\mathcal{J}$ of a normed space $(Y, \|\cdot\|_Y)$, denote by $N(\varepsilon, \mathcal{J}, \|\cdot\|_Y)$ the *minimal covering number*, i.e., the minimal number of closed balls of radius $\varepsilon$, $0 < \varepsilon < \infty$, (w.r.t. $\|\cdot\|_Y$) needed to cover $\mathcal{J}$. In accordance, let $H(\varepsilon, \mathcal{J}, \|\cdot\|_Y) = \log N(\varepsilon, \mathcal{J}, \|\cdot\|_Y)$ be the *metric entropy* of the set $\mathcal{J}$, where log denotes the natural logarithm.

The following lemma gives a relationship between metric entropy and entropy numbers:

**Lemma 10.** *Let $0 < \alpha < \infty$ and let $\mathcal{J}$ be a totally bounded (non-empty) subset of a normed space $(Y, \|\cdot\|_Y)$ satisfying*

$$e(k, \mathcal{J}, \|\cdot\|_Y) \sim k^{-1/\alpha}.$$

*We then have for the metric entropy*

$$H(\varepsilon, \mathcal{J}, \|\cdot\|_Y) \sim \varepsilon^{-\alpha}.$$

*Proof.* The inequality $H(\varepsilon) \leq C_1 \varepsilon^{-\alpha}$ is part of the proof of Theorem 1 in [12]. The lower bound follows from an obvious inversion of the argument. □

We next state a special case of more general results due to Haroske and Triebel [10]. Here we use weighted Besov spaces $B^s_{pq}(\mathbb{R}^d, \langle x \rangle^{-\gamma})$ defined in Section 4.2 in [6], see also Definition 2 in [12]. Note that $B^s_{pq}(\mathbb{R}^d, \langle x \rangle^{-0}) = B^s_{pq}(\mathbb{R}^d)$.

**Proposition 11** (Haroske and Triebel). *Suppose $p, q_1, q_2 \in [1, \infty]$, $s - d/p + d/2 > 0$. Then $B^s_{pq_1}(\mathbb{R}^d)$ is imbedded into $B^0_{2q_2}(\mathbb{R}^d, \langle x \rangle^{-\gamma})$ for every $\gamma \geq 0$. If $\gamma > 0$, the imbedding is even compact, in which case the entropy numbers of this imbedding satisfy*

$$e\left(k, id(U_{B^s_{pq_1}(\mathbb{R}^d)}), \left\|(\cdot)\langle x\rangle^{-\gamma}\right\|_{0,2,q_2,\lambda}\right) \sim k^{-1/\alpha}$$

*for all $k \in \mathbb{N}$ where $\alpha = d/s$ if $\gamma > s - d/p + d/2$ and $\alpha = (\gamma/d + 1/p - 1/2)^{-1}$ if $\gamma < s - d/p + d/2$.*

*Proof.* The first imbedding follows from Theorem 4.2.3 in [6]. The remaining claims of the proposition are proved in Theorem 4.1 in [10] for complex Besov spaces noting that the norms used in that reference are equivalent to the weighted norm $\|(\cdot)\langle x\rangle^{-\gamma}\|_{0,2,q_2,\lambda}$ used here; cf. Theorem 4.2.2 in [6]. The proposition for real Besov spaces follows from Lemma 1 in [12], see also the proof of Proposition 2 in the latter paper. □

Finally, we obtain the following corollary. [Here, and in other proofs of the paper, we use the obvious fact that metric entropy is not increased under Lipschitz-transformations between normed spaces (e.g., linear and continuous mappings); cf. also Lemma 2 in [12]].

**Corollary 12.** *Suppose $p, q \in [1, \infty]$, $s - d/p + d/2 > 0$ and $\gamma > 0$. We then have that*

(2) $$H(\varepsilon, U_{B^s_{pq}(\mathbb{R}^d)}, \left\|(\cdot)\langle x\rangle^{-\gamma}\right\|_{2,\lambda}) \sim \varepsilon^{-\alpha}$$

*where $\alpha = d/s$ if $\gamma > s - d/p + d/2$ and $\alpha = (\gamma/d + 1/p - 1/2)^{-1}$ if $\gamma < s - d/p + d/2$.*

*Proof.* We have

(3) $$H(\varepsilon, U_{B^s_{pq}(\mathbb{R}^d)}, \left\| (\cdot) \langle x \rangle^{-\gamma} \right\|_{0,2,q',\lambda}) \sim \varepsilon^{-\alpha}$$

for every $1 \leq q' \leq \infty$ by Proposition 11 and Lemma 1 above. Since

(4) $$\|f\|_{0,2,\infty,\lambda} \lesssim \|f\|_{2,\lambda} \lesssim \|f\|_{0,2,1,\lambda}$$

holds for all $f \in [f] \in B^0_{21}(\mathbb{R}^d) \supseteq B^s_{pq}(\mathbb{R}^d)$ by 2.5.7/1 in [17], we have (2) by using (3) to construct upper ($q' = 1$) and lower ($q' = \infty$) bounds for $H(\varepsilon, U_{B^s_{pq}(\mathbb{R}^d)}, \|(\cdot) \langle x \rangle^{-\gamma}\|_{2,\lambda})$. □

## Acknowledgement

The author wishes to thank Evarist Giné and Benedikt M. Pötscher for very helpful discussions (in particular about Theorem 7) and comments on a preliminary version of the paper.